\documentclass[reqno,12pt,draft]{amsart}

\usepackage{amsmath}
\usepackage{amscd}
\usepackage{amssymb}
\usepackage{enumerate}
\usepackage{amsfonts}
\usepackage{graphicx}
\usepackage[all]{xypic}
\usepackage{mathrsfs}

\DeclareFontEncoding{OT2}{}{}

\numberwithin{equation}{section}

\newcommand{\rarab}[1]{\stackrel{#1}{\longrightarrow}}
\newcommand{\rar}{\longrightarrow}

\newcommand{\al}{\alpha}  \newcommand{\ga}{\gamma}
  
\newcommand{\la}{\lambda} 
 
\newcommand{\eps}{\epsilon} 
  
 \newcommand{\Om}{\Omega}

\newcommand{\bC}{{\mathbb C}}

\newcommand{\bF}{{\mathbb F}}

\newcommand{\bN}{{\mathbb N}}

\newcommand{\bQ}{{\mathbb Q}}
\newcommand{\bR}{{\mathbb R}}

\newcommand{\bZ}{{\mathbb Z}}

\newcommand{\cF}{{\mathcal F}}

\newcommand{\cK}{{\mathcal K}}

\newcommand{\fa}{{\mathfrak a}}

\newcommand{\fg}{{\mathfrak g}}
\newcommand{\fh}{{\mathfrak h}}

\newcommand{\fk}{{\mathfrak k}}

\newcommand{\fz}{{\mathfrak z}}

\newcommand{\sU}{{\mathscr U}}
\newcommand{\sV}{{\mathscr V}}

\newcommand{\Gh}{\widehat{G}}

\newcommand{\Kh}{\widehat{K}}

\newcommand{\abs}[1]{\vert #1\vert}

\newcommand{\Ker}{\operatorname{Ker}}

\newcommand{\Hom}{\operatorname{Hom}}

\newcommand{\tr}{\operatorname{tr}}
\newcommand{\ad}{\operatorname{ad}}
\newcommand{\Ad}{\operatorname{Ad}}

\newcommand{\Supp}{\operatorname{supp}}
\newcommand{\supp}{\operatorname{supp}}

\newcommand{\Fun}{\operatorname{Fun}}

\newcommand{\Ind}{\operatorname{Ind}}

\newcommand{\res}{\operatorname{res}}

\newcommand{\tens}{\otimes}

\newcommand{\st}{\,\big\vert\,}

\newcommand{\mbr}{\medbreak}

\newcommand{\bra}{\langle}
\newcommand{\ket}{\rangle}

\newcommand{\pair}[1]{\bra #1\ket}

\newcommand{\cst}{\bC^\times}

\newtheorem{thm}{Theorem}[section]

\newtheorem{lem}[thm]{Lemma}
\newtheorem{prop}[thm]{Proposition}

\theoremstyle{remark}

\newtheorem{rems}[thm]{Remarks}

\newtheorem{defin}[thm]{Definition}

\newcommand{\bbe}{\begin{equation}}
\newcommand{\ee}{\end{equation}}

\newcommand{\Lie}{\operatorname{Lie}}

\newcommand{\xt}{\widetilde{x}}
\newcommand{\yt}{\widetilde{y}}

\title[Orbit method and Brown's theorem]{The orbit
method for profinite groups and \\ a $p$-adic analogue of Brown's
theorem}

\author[Mitya Boyarchenko and Maria Sabitova]{Mitya
Boyarchenko\address{\hskip-1.15em Mitya Boyarchenko: Department of
Mathematics, \hfill\newline University of Chicago, Chicago, IL
60637 \hfill\newline e-mail: mitya@math.uchicago.edu} \and Maria
Sabitova
\address{\hskip-1.15em Maria Sabitova: Department of Mathematics, \hfill\newline
University of Illinois at Urbana-Champaign, Urbana, IL 61801 \hfill\newline
e-mail: sabitova@math.uiuc.edu} }

\date{July 30, 2006}

\thanks{The research of M.B. was partially supported by NSF grant
DMS-0401164.}

\begin{document}

\begin{abstract}
We develop an approach to the character theory of certain classes of finite and
profinite groups based on the construction of a Lie algebra associated to such
a group, but without making use of the notion of a polarization which is
central to the classical orbit method. Instead, Kirillov's character formula
becomes the fundamental object of study. Our results are then used to produce
an alternate proof of the orbit method classification of complex irreducible
representations of $p$-groups of nilpotence class $<p$, where $p$ is a prime,
and of continuous complex irreducible representations of uniformly powerful
pro-$p$-groups (with a certain modification for $p=2$). As a main application,
we give a quick and transparent proof of the $p$-adic analogue of Brown's
theorem, stating that for a nilpotent Lie group over $\bQ_p$ the Fell topology
on the set of isomorphism classes of its irreducible representations coincides
with the quotient topology on the set of its coadjoint orbits.
\end{abstract}

\maketitle




\section*{Introduction}

The orbit method was originally discovered in the late 1950s --
early 1960s by Alexandre Kirillov \cite{kirillov} for connected
and simply connected nilpotent Lie groups. If $G$ is such a group
and $\fg$ is its Lie algebra, this method provides an explicit
bijection between the {\em unitary dual} $\Gh$ of $G$, i.e., the
set of equivalence classes of unitary irreducible representations
of $G$, and the set $\fg^*/G$ of orbits of the induced action of
$G$ on $\fg^*$ (called {\em coadjoint orbits}). A major ingredient
of this theory is {\em Kirillov's character formula}. Roughly
speaking, it states that if $\Om\subset\fg^*$ is a coadjoint orbit
and $\rho_\Om\in\Gh$ is the corresponding representation, then the
character of $\rho_\Om$, viewed as a generalized function on $G$,
is the pullback via the logarithm map $\log:G\rar{}\fg$ of the
inverse Fourier transform of a suitably normalized $G$-invariant
measure on $\fg^*$ supported on $\Om$.

\mbr

Since then Kirillov's approach has been extended to many other classes of
groups: nilpotent $p$-adic Lie groups \cite{moore}, $p$-groups of nilpotence
class $<p$ (beginning with \cite{kazhdan}), and uniformly powerful (or uniform,
for short) pro-$p$-groups \cite{howe,jaikin-zapirain-zeta}, to name the ones
that will appear in this paper\footnote{The orbit method can also be applied,
with suitable changes, to solvable Lie groups; moreover, its philosophy extends
essentially to all Lie groups, and even beyond them. However, these
generalizations lie in a different direction from the ones considered in this
article.}. Each such extension usually involves two modifications: one has to
work with a correct analogue of a ``unitary irreducible representation'' in
each context, and one has to find an appropriate version of the Lie algebra
construction. For example, if $G$ is a $p$-adic Lie group, its Lie algebra is
defined as usual, but $\Gh$ has to be understood as the set of isomorphism
classes of irreducible complex ``algebraic'' (or, in a different terminology,
``smooth'') representations of $G$. On the other hand, if $G$ is a $p$-group of
nilpotence class $<p$ (respectively, a uniform pro-$p$-group), then $\Gh$ has
to be understood as the set of isomorphism classes of continuous complex
irreducible representations of $G$, and the usual Lie algebra construction is
replaced by a construction of Lazard which produces a finite Lie ring \cite{kh}
(respectively, a uniform Lie algebra over $\bZ_p$ \cite{pro-p}) associated to
$G$.

\mbr

After these modifications have been made, the theory follows the
pattern of Kirillov's original approach (modulo various technical
difficulties). Namely, in each case the underlying additive group
of $\fg$ has a natural topology, and $\fg^*$ can be identified
with the {\em Pontryagin dual} of $\fg$. Given an element
$f\in\fg^*$, one looks for a \emph{polarization} of $\fg$ at $f$,
i.e., a \emph{Lie subalgebra} $\fh\subseteq\fg$ which has the
property that $f$ is trivial on $[\fh,\fh]$, and which is maximal
among all \emph{additive subgroups} of $\fg$ with this property.
Polarizations always exist, and if $H$ is the subgroup of $G$
corresponding to a polarization $\fh$, then $f$ induces a
$1$-dimensional character $\chi_f$ of $H$ and we can form the
induced representation $\rho_f=\Ind_H^G\chi_f$. The theorem is
that this representation is always irreducible; its isomorphism
class depends only on the $G$-orbit of $f$; and, finally, every
$\rho\in\Gh$ arises in this way from a unique $G$-orbit
$\Om\subset\fg^*$. This description of $\rho_f$ is then used to
prove Kirillov's character formula (or a suitable analogue
thereof).

\mbr

In reality, one needs to be more careful with uniformly powerful
pro-$p$-groups when $p=2$. The problem that arises here is that an
element $f\in\fg^*$ which is trivial on $[\fg,\fg]$ may not induce
a $1$-dimensional character of G. (We thank A.~Jaikin-Zapirain for
explaining this to us.) Thus in this case the approach to the
orbit method has to be somewhat modified
(cf.~\cite{jaikin-zapirain-zeta}); however, the basic idea remains
the same.

\mbr

An important feature of all four situations mentioned above is
that both $\Gh$ and $\fg^*/G$ are equipped with a natural
topology. The topology on the former is the so-called {\em Fell
topology} (see \S\ref{ss:fell}). The topology on the latter is the
quotient of the standard (compact-open) topology on $\fg^*$.
Moreover, in all four cases the orbit method bijection turns out
to be a homeomorphism. This is a nontrivial result which has
useful applications. For an interesting application in the
$p$-adic setting we refer the reader to \cite{gelfand-kazhdan}. In
the setting of real Lie groups this statement was originally
conjectured by Kirillov, who also proved that the bijection
$\fg^*/G\rar{}\Gh$ is continuous. The proof that this bijection is
also open is substantially more difficult, and was given by Ian
Brown about 10 years later in \cite{brown}. While it may be
possible to adapt Brown's argument to a $p$-adic nilpotent Lie
group $G$ (to the best of our knowledge, this has never been
done), we present in Section \ref{s:brown} a completely different
proof (following a suggestion of V.~Drinfeld) which is based on
the fact that $G$ is an increasing union of a sequence of open
uniform pro-$p$-subgroups (see Lemma \ref{l:union}); this is the
main new result of our paper. Our proof seems to be much shorter
and more transparent than Brown's proof, and we hope that it is
easier to understand. On the other hand, it is not clear to us
whether this approach has an analogue for real Lie groups.

\mbr

Another new result in our paper is a theorem we call the
``abstract orbit method''. It arose from an approach to the orbit
method for finite nilpotent groups (of sufficiently small
nilpotence class) that we also learned from V.~Drinfeld. It was
natural to try to see if this approach can be extended to uniform
pro-$p$-groups, and, more generally, to find the minimal set of
assumptions under which this method can be used. The answer is
given in Section \ref{s:abstract}, and in Section \ref{s:pro-p} we
show that our ``abstract orbit method'' can indeed be used to
classify (continuous) complex irreducible representations of
$p$-groups of nilpotence class $<p$ and of uniform pro-$p$-groups
(with a certain modification for $p=2$). The main difference with
the classical approach is that we never mention polarizations. In
particular, in the abstract setting one does not even need a Lie
bracket on $\fg$. Instead, we prove directly that a suitable
analogue of Kirillov's character formula produces a collection of
functions on the group, parameterized by the coadjoint orbits,
which turn out to be precisely the irreducible characters of the
group.

\mbr

This approach has its advantages and disadvantages. The main
disadvantage is that, unlike the classical one, our method of
constructing irreducible characters cannot be ``upgraded'' to
yield a construction of irreducible representations. On the other
hand, it appears to be more straightforward, since one always
works directly with irreducible characters, whereas the motivation
behind the notion of a polarization comes from areas of
mathematics outside of representation theory. However, a much more
significant advantage is that the method explained in our paper
has an analogue in the geometric representation theory for
unipotent groups, whereas the classical method does not have such
an analogue (at least not in any obvious sense), for in the
geometric setting polarizations cease to exist in general. A
proper discussion of this remark is beyond the scope of our paper,
and instead we refer the reader to \cite{intro}.

\subsection*{Acknowledgements} This paper owes its existence to
lectures of Vladimir Drinfeld and our private discussions with
him. In particular, he explained to us the approach to the orbit
method for finite nilpotent groups which gave rise to our
``abstract orbit method'' theorem. He also motivated our main
result by asking if an analogue of Brown's theorem for $p$-adic
nilpotent Lie groups can be proved using the orbit method for
uniformly powerful pro-$p$-groups.

\mbr

In addition, we would like to thank Michael Geline and George
Glauberman who told us about Lazard's construction in the setting
of $p$-groups and uniformly powerful pro-$p$-groups, and suggested
references \cite{kh} and \cite{pro-p}, respectively. We are also
grateful to Ben Wieland for referring us to \cite{howe}, to Adam
Logan for drawing our attention to \cite{jaikin-zapirain-zeta},
and especially to Andrei Jaikin-Zapirain for helpful e-mail
correspondence and for pointing out a mistake in an earlier
version of our paper.


\section{Abstract orbit method}\label{s:abstract}

\subsection{The statement}\label{ss:abstract-statement}
For every profinite group $\Pi$ we denote by $\mu_\Pi$ the unique
Haar measure on $\Pi$ such that $\mu_\Pi(\Pi)=1$. We define the
\emph{convolution} of two complex-valued $L^{2}$-functions $f_1$
and $f_2$ on $\Pi$ by the formula
\[
(f_1*f_2)(\ga)=\int_G f_1(h)f_2(h^{-1}\ga)\,d\mu_\Pi(h),\qquad
\ga\in \Pi.
\]
We write $\Fun(\Pi)$ for the space of complex-valued functions on
$\Pi$ that are bi-invariant with respect to a sufficiently small
open subgroup of $\Pi$. It is clear that $\Fun(\Pi)\subseteq
L^2(\Pi)$ is closed under convolution, which makes $\Fun(\Pi)$ an
associative $\bC$-algebra (it is unital if and only if $\Pi$ is
finite, and commutative if and only if $\Pi$ is commutative). The
subspace $\Fun(\Pi)^\Pi\subseteq\Fun(\Pi)$ of $\Pi$-invariant
functions, where $\Pi$ acts on itself by conjugation, is also
closed under convolution (see Lemma
\ref{l:characters-idempotents}), and in fact coincides with the
\emph{center} of $\Fun(\Pi)$; in particular, $\Fun(\Pi)^\Pi$ is
always commutative.
\begin{thm}[Abstract orbit method]\label{t:abstract}
Let $G$ be a profinite group, and suppose that there exist an
abelian profinite group $\fg$ and a homeomorphism
$\exp:\fg\rar{}G$ such that the following two conditions hold:
\begin{enumerate}[$($i$)$]
\item for each $g\in G$, the map $\Ad g:\fg\rar{}\fg$ given by
$x\longmapsto\log(g\exp(x)g^{-1})$ is a group automorphism, where
we write $\log$ for $\exp^{-1}$; and
\item the pullback map $\exp^*:\Fun(G)^G\rarab{\simeq}\Fun(\fg)^G$
commutes with convolution.
\end{enumerate}
Then each $G$-orbit $\Om\subset\fg^*$ is finite, and there is a
bijection between $\fg^*/G$ and $\Gh$ such that the irreducible
character $\chi$ of $G$ corresponding to an orbit
$\Om\subset\fg^*$ is given by
\begin{equation}\label{e:charform-profinite}
\chi(e^x)=\abs{\Om}^{-1/2} \sum_{f\in\Om} f(x).
\end{equation}
\end{thm}
Here, as in the introduction, $\fg^*$ denotes the Pontryagin dual
of $\fg$, which, since $\fg$ is compact, coincides with the group
of continuous homomorphisms of $\fg$ into $\cst$, and has the
discrete topology. The action of $G$ on $\fg^*$ is induced by its
action on $\fg$ via $\Ad$. Note that every finite group can be
viewed as a profinite one (with the discrete topology), so our
definitions and the theorem are valid for finite $G$ and $\fg$ as
well.
\begin{rems}
\begin{enumerate}[(1)]
\item In practice, if one wants to apply Theorem \ref{t:abstract}
to a specific group $G$, the main difficulty lies in verifying
assumption (ii), as we will see in Section \ref{s:pro-p}.
\item As we have already noted in the introduction, one should
observe that $\fg$ is not required to have a Lie bracket in the
statement of the theorem. Unfortunately, we do not know of any
example where the assumption of the theorem is satisfied for some
profinite group $G$, but $\fg$ does not arise from some sort of a
Lie algebra construction. It would be very interesting to find
such an example.
\item Formula \eqref{e:charform-profinite} implies that
$\abs{\Om}^{1/2}=\chi(1)$ is an integer for every $\Om\in\fg^*/G$,
i.e., the order of every coadjoint orbit is a full square. In the
generality of the theorem, this is the only proof of this fact
known to us.
\end{enumerate}
\end{rems}

\subsection{Auxiliary results}\label{ss:abstract-auxiliary} Until the
end of the section we fix $G$, $\fg$ and $\exp$ satisfying the
assumptions of the theorem. To simplify notation, we write
$X_G=\Fun(G)^G$. Thus $X_G$ is the set of all complex-valued {\em
class functions} $f$ on $G$ such that there exists a normal open
subgroup $K$ of $G$ (depending on $f$) satisfying:
$$
\forall g\in G,\quad \forall k\in K: \quad f(gk)=f(kg)=f(g).
$$
Similarly, we will write $X_\fg=\Fun(\fg)^G$, and we write
$gxg^{-1}$ in place of $(\Ad g)(x)$.
\begin{lem}\label{l:characters-idempotents}
We have $X_G\subseteq L^2(G)$, and $X_G$ is an algebra with
respect to convolution. Also, if $\chi$ is the character of a
continuous complex irreducible representation $\rho$ of $G$, then
$(\dim\rho)\cdot\chi$ is an indecomposable idempotent of $X_G$,
and every indecomposable idempotent of $X_G$ has this form.
\end{lem}

\mbr

Let us recall that an {\em indecomposable idempotent} of a
commutative ring $A$ is a nonzero idempotent $e\in A$ (i.e.,
$e\neq 0$ and $e^2=e$) which cannot be written as $e=e_1+e_2$ for
nonzero idempotents $e_1,e_2\in A$ satisfying $e_1\cdot e_2=0$.

\begin{proof}[Proof of Lemma $\ref{l:characters-idempotents}$] For any $f_1,f_2\in X_G$
let $K$ be a normal open subgroup of $G$ such that both $f_1$ and
$f_2$ are constant on the cosets of $K$ in $G$. (Clearly, such a
$K$ exists.) Then $f_1$ and $f_2$ can be considered as class
functions on $G/K$ (denoted respectively as $\bar{f_1}$ and
$\bar{f_2}$), and for each $g\in G$, we have
\[
(f_1*f_2)(g)=(\bar{f_1}*\bar{f_2})(\bar{g}),
\]
where $\bar{g}$ denotes the image of $g$ in $G/K$. On the other
hand, it is well known that every complex irreducible
representation of $G$ is finite dimensional (because $G$ is
compact), and hence has finite kernel (because $G$ is totally
disconnected). This implies that it is enough to prove the lemma
for a finite group $G$.

\mbr

Let $G$ be finite. Then clearly the set $\Fun(G)$ of all functions
on $G$ is an algebra with convolution as multiplication and
$X_G\subseteq \Fun(G)$ is the set of all class functions on $G$.
Define the map $\psi:\bC[G]\rar \Fun(G)$ via
\begin{displaymath}
\psi(g)(h)= \left\{\begin{array}{cc}
0, & g\ne h \\
\abs{G}, & g=h
\end{array}\right. \quad g,h\in G,
\end{displaymath}
and for any $x=\sum_{g\in G}n_gg$ $(n_g\in\bC)$,
$$
\psi(x)=\sum_{g\in G}\overline{n_g}\psi(g),
$$
where $\overline{n_g}$ denotes the complex conjugate of $n_g$. It
is easy to see that $\psi$ is a ring isomorphism and that the
inverse image of $X_G$ in $\bC[G]$ is center $ZG$ of $\bC[G]$.
This implies that $X_G\subseteq \Fun(G)$ is a subalgebra and that
$e\in ZG$ is an indecomposable idempotent if and only if $\psi(e)$
is one.

\mbr

Recall that there is a one-to-one correspondence between
indecomposable idempotents of $ZG$ and irreducible representations
of $G$, such that for every indecomposable idempotent $e\in ZG$,
the corresponding irreducible representation $\rho_{e}$ of $G$ has
the property that the left regular representation of $G$ on
$e\cdot\bC[G]$ is isomorphic to a multiple of $\rho_e$. Thus, it
is enough to show that if $e\in ZG$ is an indecomposable
idempotent corresponding to an irreducible representation $\rho_e$
of $G$, then $(\dim\rho_e)^{-1}\cdot\psi (e)$ is the character of
$\rho_e$. Moreover, since $\psi(e)$ is a class function, it is
enough to show that for any two indecomposable idempotents
$e,\tilde{e}\in ZG$ we have \bbe\label{eq:1} \langle
\psi(e),\chi_{\tilde{e}}\rangle= \left\{\begin{array}{cc}
0, & e\ne \tilde{e} \\
\dim\rho_e, & e=\tilde{e}
\end{array}\right.
\ee where $\langle \cdot\, ,\cdot\rangle$ is the inner product on
$L^2(G)=\Fun(G)$ and $\chi_{\tilde{e}}$ is the character of
$\rho_{\tilde{e}}$.

\mbr

Let $e=\sum_{g\in G}n_gg$, $n_g\in\bC$. Then \bbe\label{eq:2}
\langle
\psi(e),\chi_{\tilde{e}}\rangle=\frac{1}{\abs{G}}\cdot\sum_{g\in
G}\overline{\psi(e)(g)}\cdot\chi_{\tilde{e}}(g)= \sum_{g\in
G}n_g\chi_{\tilde{e}}(g)=\chi_{\tilde{e}}(e), \ee where by
extending $\chi_{\tilde{e}}$ by linearity we consider
$\chi_{\tilde{e}}$ as a function on $\bC[G]$. Now let
$V_{\tilde{e}}$ be a representation space for $\rho_{\tilde{e}}$.
Then, as was mentioned above, for some $n\in\bN$ there exists a
$\bC[G]$-module isomorphism
$$
\tilde{e}\cdot\bC[G]\cong V_{\tilde{e}}^n,
$$
hence
\begin{displaymath}
\chi_{\tilde{e}}(e)= \left\{\begin{array}{cc}
0, & e\ne \tilde{e} \\
\dim\rho_e, & e=\tilde{e}
\end{array}\right.
\end{displaymath}
which together with \eqref{eq:2} gives \eqref{eq:1}.
\end{proof}

\begin{lem}\label{l:coad-finite}
Every $G$-orbit in $\fg^*$ is finite.
\end{lem}
\begin{proof}
Fix $f\in\fg^*$. Since $\fg$ is profinite and $f:\fg\to\cst$ is
continuous, it follows that $f$ has finite image, hence $\Ker f$
is open. Thus there exists an open normal subgroup $K\subseteq G$
such that $\log(K)\subseteq\Ker f$. Hence every element in the
$G$-orbit of $f$ vanishes on $\log(K)$. But the set of elements of
$\fg^*$ that vanish on $\log(K)$ is finite, since there exists an
open additive subgroup $\fa\subseteq\fg$ such that
$\fa\subseteq\log(K)$, and $(\fg/\fa)^*$ is finite.
\end{proof}

\begin{lem}\label{l:characters-orbits}
Write $\mu=\mu_\fg$. For each $\Om\in\fg^*/G$ define the following
function on $\fg$:
$$
\chi_{\Om}(x)=\frac{1}{\abs{\Om}^{1/2}}\sum_{f\in\Om}f(x),\quad
x\in\fg.
$$
Then $\chi_\Om\in X_\fg$, and for any two $\Om,\Om'\in\fg^*/G$ we
have
\begin{displaymath}
\pair{\chi_{\Om},\chi_{\Om'}}= \left\{\begin{array}{cc}
1, & \Om=\Om' \\
 0, & \Om\ne  \Om'
\end{array}\right.
\end{displaymath}
where $\pair{\chi_{\Om},\chi_{\Om'}}=
\int_{\fg}\chi_{\Om}(x)\overline{\chi_{\Om'}(x)}d\mu(x)$.
\end{lem}
\begin{proof}
It is clear that $\chi_\Om\in X_\fg$. Moreover, the functions $f\in\fg^*$ are
well known to be orthonormal with respect to the inner product
$\pair{\cdot\,,\cdot}$. Indeed, this is simply the orthogonality of irreducible
characters for the compact abelian group $\fg$. This implies in particular that
the $L^2$ norm of the function $\sum_{f\in\Om} f$ is equal to
$\abs{\Om}^{1/2}$. Since distinct orbits are disjoint, the statement of the
lemma follows immediately.
\end{proof}

\subsection{Proof of Theorem
\ref{t:abstract}}\label{ss:abstract-proof} Let us write
$L^2(\fg^*)$ for the space of square-summable functions on $\fg^*$
(where $\fg^*$ is equipped with the counting measure). Recall that
the Fourier transform provides an isomorphism
\[
\cF:L^2(\fg)\rarab{\simeq} L^2(\fg^*),
\]
\[
(\cF f)(\phi)=\int_\fg f(h)\overline{\phi(h)}d\mu(h),\quad f\in
L^2(\fg),\ \phi\in \fg^*,
\]
which intertwines convolution with pointwise multiplication
(whenever the two operations are defined). Let $\chi_\Om$ be the
functions on $\fg$ defined in Lemma \ref{l:characters-orbits}.
Given $\Om\in\fg^*/G$, we now show that
$\cF(\abs{\Om}^{1/2}\cdot\chi_{\Om})$ is the characteristic
function of $\Om$. Let $\Om'\in\fg^*/G$ be a second orbit
(possibly the same as $\Om$), and let $\phi\in\Om'$. Then we have
$$
\cF(\chi_{\Om})(\phi)=\int_{\fg}\chi_{\Om}(x)\overline{\phi(x)}d\mu(x).
$$
This implies that
\begin{eqnarray*}
&&
\pair{\chi_{\Om},\chi_{\Om'}}=\int_{\fg}\chi_{\Om}(x)\overline{\chi_{\Om'}(x)}d\mu(x)=\\
&& \frac{1}{\abs{\Om'}^{1/2}\abs{G^{\phi}}}\sum_{g\in
G}\int_{\fg}\chi_{\Om}(x)\overline{\phi(g^{-1}xg)}d\mu(x)=\\
&&
\frac{\abs{G}}{\abs{\Om'}^{1/2}\abs{G^{\phi}}}\int_{\fg}\chi_{\Om}(x)\overline{\phi(x)}d\mu(x)=
\abs{\Om'}^{1/2}\cdot\cF(\chi_{\Om})(\phi).
\end{eqnarray*}
Thus, by Lemma \ref{l:characters-orbits},
$\cF(\abs{\Om}^{1/2}\cdot\chi_{\Om})$ is the characteristic
function of $\Om$, and is therefore an indecomposable idempotent
in the algebra of $G$-invariant functions on $\fg^*$ with respect
to pointwise multiplication. Hence
$\abs{\Om}^{1/2}\cdot\chi_{\Om}$ is an indecomposable idempotent
of the algebra $X_G$ (with respect to convolution), which together
with Lemma \ref{l:characters-idempotents} proves the
correspondence between irreducible representations of $G$ and
coadjoint orbits of $\fg^*.$ Furthermore, if $\chi$ is the
irreducible character of $G$ corresponding to $\Om$, then
$$
\abs{\Om}^{1/2}\cdot\chi_{\Om}=\dim\rho\cdot\exp^*(\chi)
$$
by Lemma \ref{l:characters-idempotents}, hence by evaluating these
functions at 0 we see that $\abs{\Om}^{1/2}=\dim\rho$, and
consequently $\chi_{\Om}=\exp^*(\chi)$. This completes the proof
of Theorem \ref{t:abstract}.


\section{Application to (pro-)$p$-groups}\label{s:pro-p}

\subsection{Notation and terminology} Until the end of the paper,
$p$ will denote a fixed prime number. Even though there exists a
group-theoretic definition of a uniformly powerful (or, for
brevity, ``uniform'') pro-$p$-group \cite{pro-p}, for our purposes
it is more convenient to use the Lie-theoretic definition, which
is also more transparent. We define a \emph{uniform Lie algebra}
to be a Lie algebra over the ring $\bZ_p$ of $p$-adic integers
which is free of finite rank as a $\bZ_p$-module and satisfies
$[\fg,\fg]\subseteq p\cdot\fg$ (respectively, $[\fg,\fg]\subseteq
4\cdot\fg$ when $p=2$). Given a uniform Lie algebra $\fg$, we
equip it with the topology induced by the standard topology on
$\bZ_p$, and we define a topological group $G:=\exp\fg$ to be the
underlying topological space of $\fg$ equipped with a group
operation given by the \emph{Campbell-Hausdorff series}
\begin{equation}\label{e:C-H}
CH(x,y) = \log\bigl(\exp(x)\exp(y)\bigr) = \sum_{i=1}^\infty
CH_i(x,y).
\end{equation}
\begin{rems}
\begin{enumerate}[(1)]
\item A priori, $CH(x,y)$ is viewed as an element of
$\bQ\langle\langle x,y\rangle\rangle$, the algebra of formal
noncommutative power series in the variables $x$ and $y$ with
coefficients in $\bQ$, and $CH_i(x,y)$ denotes its homogeneous
component of (total) degree $i$. However, it is well known that
$CH(x,y)$ is in fact a \emph{formal Lie series}, which means that
each term $CH_i(x,y)$ lies in the Lie subalgebra of
$\bQ\langle\langle x,y\rangle\rangle$ generated by $x$ and $y$.
\item Since the coefficients of $CH(x,y)$ involve positive powers
of $p$ in the denominator, it is not immediately obvious that
$CH(x,y)$ can even be evaluated in a uniform Lie algebra $\fg$
term-by-term. However, Michel Lazard proved, cf.~\cite{pro-p},
that the condition imposed on $\fg$ guarantees that for each
$x,y\in\fg$, we have
$CH_i(x,y)\in\fg\subseteq\bQ_p\tens_{\bZ_p}\fg$, and, in addition,
the series $CH(x,y)$ converges uniformly on $\fg$ and makes it a
topological group. This result justifies our construction of
$G=\exp\fg$.
\end{enumerate}
\end{rems}

\begin{defin}
A \emph{uniform pro-$p$-group} is a profinite group $G$ which is
isomorphic to $\exp\fg$ for some uniform Lie algebra $\fg$. If
$G\cong\exp\fg$ is such a group, we will fix an isomorphism
$\exp\fg\rarab{\simeq}G$ and denote the underlying map of sets by
$\exp:\fg\rar{}G$. By abuse of notation, we will also write
$G=\exp\fg$ and $\fg=\Lie(G)$.
\end{defin}

\mbr

Suppose now that $\fg$ is a finite Lie ring (i.e., a Lie algebra
over $\bZ$) whose order is a power of $p$, and such that $\fg$ is
nilpotent of nilpotence class $<p$. (This means that any iterated
commutator of length $\geq p$ vanishes in $\fg$.) In this case it
is rather easy to check that the Lie series $CH(x,y)$ can be
evaluated in $\fg$ term-by-term, and makes $\fg$ a $p$-group
(there is no issue of convergence of $CH$ in this setting). Again,
we denote this $p$-group by $\exp\fg$. Michel Lazard proved,
cf.~\cite{kh}, that every $p$-group $G$ of nilpotence class $<p$
arises in this way from a unique $\fg$. In this situation we will
also write $\fg=\Lie(G)$, and the underlying map of a fixed
isomorphism $\exp\fg\rarab{\simeq}G$ will be denoted by
$\exp:\fg\rar{}G$, just as for uniform pro-$p$-groups. We will
write $\log=\exp^{-1}$.

\subsection{Auxiliary results on formal Lie series}\label{ss:formal-Lie}
Throughout the rest of the section $G$ will denote a $p$-group of
nilpotence class $<p$ or a uniform pro-$p$-group, and
$\fg=\Lie(G)$ its Lie algebra. In both cases we have the
exponential map $\exp:\fg\rar{}G$, which is a homeomorphism, and
assumption (i) of Theorem \ref{t:abstract} is satisfied in this
situation. In fact, in both cases it is known that Lazard's
construction is functorial; in particular, a continuous
set-theoretic bijection $\fg\rar{}\fg$ is a Lie algebra
automorphism if and only if it is an automorphism of the group
$\exp\fg$. Another fact that will often be used implicitly in what
follows is that if $x\in\fg$ and $\ad x:\fg\rar{}\fg$ denotes the
additive map $y\mapsto [x,y]$, then $e^{\ad x}=\Ad(e^x)$ as
automorphisms of $\fg$. In order to classify the continuous
complex irreducible representations of $G$ we would like to show
that assumption (ii) of Theorem \ref{t:abstract} holds in this
setting as well. Unfortunately, as A.~Jaikin-Zapirain pointed out
to us, this is sometimes false when $p=2$ and $G$ is a uniform
pro-$2$-group; thus this case needs to be dealt with separately
(see \S\ref{ss:pro-2-groups}).

\mbr

Ignoring this issue for the moment, let us note that the main
problem with verifying assumption (ii) arises from the fact that
the convolution of functions on $\fg$ is defined using the
addition in $\fg$, while the convolution of functions on $\exp\fg$
is defined using the multiplication in $\exp\fg$, or,
equivalently, the Campbell-Hausdorff operation
$CH:\fg\times\fg\rar{}\fg$. This problem is dealt with in a very
natural way, shown to us by V.~Drinfeld: we prove that one can
write $CH(x,y)=\widetilde{x}+\widetilde{y}$, where $\widetilde{x}$
(resp., $\widetilde{y}$) is a certain Lie series in the variables
$x,y$ which is conjugate to $x$ (resp., to $y$). This formula
implies that the two convolutions of conjugation-invariant
functions on $\fg$, defined using addition and the operation $CH$,
are in fact identical, which is the content of condition (ii) of
Theorem \ref{t:abstract}. We should mention that in practice,
however, the realization of this idea for uniform pro-$p$-groups
involves certain technical difficulties.

\mbr

We now turn to precise statements. Let $E$ be a finite extension
of $\bQ_p$, let $v_p$ denote the valuation on $E$ normalized by
$v_p(p)=1$, and let $\cK\subseteq E$ be a subfield.
\begin{lem}\label{l:phi-psi}
Let $H(x,y)=\sum_{n=1}^\infty H_n(x,y)\in \cK\langle\langle
x,y\rangle\rangle$ be a formal Lie series, where $H_n(x,y)$ is
homogeneous of degree $n$, such that $H_1(x,y)=x+y$ and
$v_p(H_n)\geq -\frac{n-2}{p-1}$ for all $n\geq 2$. Then there
exist formal Lie series $\phi=\phi(x,y),\,\psi=\psi(x,y)\in
\cK\langle\langle x,y\rangle\rangle$ such that
\begin{equation}\label{e:H-phi-psi}
H(x,y) = e^{\ad\phi(x,y)}(x) + e^{\ad\psi(x,y)}(y),
\end{equation}
and if $\phi_n$, $\psi_n$ denote the degree $n$ homogeneous
components of $\phi$, $\psi$, respectively, then $v_p(\phi_n)\geq
-\frac{n-1}{p-1}$ and $v_p(\psi_n)\geq -\frac{n-1}{p-1}$ for all
$n\geq 1$.
\end{lem}
Here, by abuse of notation, we write $v_p(H_n)$ $\bigl($resp.,
$v_p(\phi_n)$ and $v_p(\psi_n)${}$\bigr)$ for the minimum among
the valuations of all coefficients of $H_n$ (resp., $\phi_n$ and
$\psi_n$).
\begin{proof}
It is easy to see that one can construct the series $\phi$ and
$\psi$ inductively. Namely, for each $n\geq 0$ let us compare the
homogeneous components of degree $n+1$ on both sides of
\eqref{e:H-phi-psi}. For $n=0$ there is nothing to check, thanks
to the assumption that $H_1(x,y)=x+y$. For each $n\geq 1$, we may
assume that all $\phi_j,\,\psi_j$ with $j<n$ have already been
found and satisfy $v_p(\phi_j),\,v_p(\psi_j)\geq
-\frac{j-1}{p-1}$. In order to find $\phi_n$ and $\psi_n$ we have
to solve an equation of the form
\begin{equation}\label{e:step-n}
[\phi_n,x] + [\psi_n,y] + \text{(something known)} = H_{n+1}(x,y),
\end{equation}
where ``something known'' is a sum of expressions of the form
\begin{equation}\label{e:known}
\frac{1}{k!}\cdot
[\phi_{j_1},[\phi_{j_2},[\dotsc[\phi_{j_k},x]\dotsc]]]
\quad\text{or}\quad \frac{1}{k!}\cdot
[\psi_{j_1},[\psi_{j_2},[\dotsc[\psi_{j_k},y]\dotsc]]]
\end{equation}
with $2\leq k\leq n$ and $j_1+j_2+\dotsb+j_k=n$. It is well known
that $v_p(k!)\leq\frac{k-1}{p-1}$ for all $k\geq 1$ (with equality
if $k$ is a power of $p$), which implies that the valuation of
each of the expressions in \eqref{e:known} is at least
\[
-\frac{k-1}{p-1}-\frac{j_1-1}{p-1} - \dotsb - \frac{j_k-1}{p-1} =
- \frac{k-1+j_1+\dotsb+j_k-k}{p-1}= - \frac{n-1}{p-1}.
\]
In addition, we have $v_p(H_{n+1})\geq -\frac{n-1}{p-1}$ by
assumption. This immediately implies that there exist homogeneous
Lie polynomials $\phi_n=\phi_n(x,y)$, $\psi_n=\psi_n(x,y)$ of
degree $n$ which solve \eqref{e:step-n} and satisfy
$v_p(\phi_n),\,v_p(\psi_n)\geq -\frac{n-1}{p-1}$, completing the
induction.
\end{proof}
This result suffices to prove the orbit method correspondence when
$p\geq 5$. To treat the case $p=3$ we need the following
variation:
\begin{lem}\label{l:ugly}
In the situation of Lemma $\ref{l:phi-psi}$, assume that $p=3$ and that
\[
v_3(H_n)\geq -\frac{6n-10}{7} \qquad\forall\,n\geq 2.
\]
Then the conclusion of Lemma $\ref{l:phi-psi}$ holds with
$v_3(\phi_n),\,v_3(\psi_n)\geq -\frac{6n-4}{7}$ for $n\geq 1$.
\end{lem}
\begin{proof}
We follow the proof of Lemma \ref{l:phi-psi} almost word-for-word;
the only step that needs to be changed is the estimation of the
valuations of the coefficients of the expressions \eqref{e:known}.
We have $v_3(k!)\leq(k-1)/2$, and therefore, by the induction
assumption, each of the valuations in question is at least
\[
-\frac{k-1}{2}-\frac{6j_1-4}{7}-\dotsb-\frac{6j_k-4}{7} =
-\frac{k-1}{2}-\frac{6n-4k}{7} = -\frac{6n-4}{7}+\frac{k-1}{14}.
\]
Since $k-1>0$, this finishes the induction in the same way as
before.
\end{proof}

\subsection{The orbit method when $p\geq 3$}\label{ss:p-not-2}
In this subsection we treat the orbit method for a group $G$ which
is either a $p$-group of nilpotence class $<p$ or a uniform
pro-$p$-group with $p\geq 3$. Since the orbit method obviously
works for commutative $2$-groups, there is no harm in assuming
that $p\geq 3$ for finite $G$ as well.
\begin{prop}\label{p:phi-psi}
Assume that $p\geq 3$, let $G$ be as above, and let $\fg=\Lie(G)$.
Then there exist formal Lie series
$\phi(x,y),\,\psi(x,y)\in\bQ\langle\langle x,y\rangle\rangle$
which can be evaluated term-by-term in $\fg$, converge uniformly
for $x,y\in\fg$ when $\fg$ is uniform, and satisfy
\begin{equation}\label{e:phi-psi-general}
\log\bigl(\exp(x)\exp(y)\bigr) =
e^{\ad\phi(x,y)}(x)+e^{\ad\psi(x,y)}(y) \qquad \forall\,x,y\in\fg.
\end{equation}
\end{prop}
\begin{proof}
Let us recall the Campbell-Hausdorff series
$CH(x,y)\in\bQ\langle\langle x,y\rangle\rangle$ defined by
\eqref{e:C-H}. The key fact about the coefficients of $CH(x,y)$
that we will need is the following result (see \cite{pro-p},
p.~123): for every prime $p$,
\begin{equation}\label{e:val-CH}
v_p(CH_n)\geq -\frac{n-1}{p-1} \qquad \forall\,n\geq 1.
\end{equation}
Suppose first that $G$ is finite, and let
$H(x,y)\in\bQ\langle\langle x,y\rangle\rangle$ denote the Lie
polynomial obtained by discarding all homogeneous components of
$CH(x,y)$ of degrees $\geq p$. Then $G$ is isomorphic to $\fg$
equipped with the operation given by $H$. Since $v_p(H_n)\in\bZ$,
it follows from \eqref{e:val-CH} that all coefficients of $H$ lie
in $\bZ_p\cap\bQ$. Thus the assumption of Lemma \ref{l:phi-psi} is
satisfied with $\cK=\bQ$. Let $\phi'$, $\psi'$ denote formal Lie
series satisfying the conclusion of the lemma, and let $\phi$,
$\psi$ be the Lie polynomials obtained from $\phi'$ and $\psi'$,
respectively, by discarding all homogeneous components of degrees
$\geq p-1$. (Note a slight change of our notation!) Since the
valuation of each coefficient of $\phi$ and $\psi$ must be an
integer, and since $-\frac{n-1}{p-1}>-1$ for $1\leq n\leq p-2$, we
see that the coefficients of $\phi$ and $\psi$ lie in
$\bZ_p\cap\bQ$. Thus $\phi$ and $\psi$ can be evaluated in $\fg$.
In addition, since $\fg$ has nilpotence class $\leq p-1$, the
conclusion of Lemma \ref{l:phi-psi} implies that
\eqref{e:phi-psi-general} holds.

\mbr

Next we assume that $G$ is a uniform pro-$p$-group. By abuse of
notation, we define a ``valuation''
$v_p:\bQ_p\tens_{\bZ_p}\fg\rar{}\bZ\cup\{\infty\}$ by
$v_p(x)=\sup\bigl\{ r\in\bZ \,\bigl\lvert\, p^{-r}x\in\fg
\bigr\}$. It is well known that a series $\sum_{n=1}^\infty x_n$
in $\fg$ converges if and only if $v_p(x_n)\to\infty$ as
$n\to\infty$.

\mbr

Let us consider the case $p=3$. Put
$H(x,y)=CH(x,y)\in\bQ\langle\langle x,y\rangle\rangle$. It follows
from \eqref{e:val-CH} that $v_3(H_2)\geq 0$ and $v_3(H_n)\geq
-\frac{n-1}{2}$ for all $n\geq 3$. Since
$-\frac{6n-10}{7}\leq-\frac{n-1}{2}$ for all $n\geq 3$, we see
that $H$ satisfies the assumption of Lemma \ref{l:ugly} with
$\cK=\bQ$. Let $\phi$, $\psi$ denote the formal Lie series
satisfying the conclusion of the lemma. Since $\fg$ is uniform, we
see that for all $x,y\in\fg$ and all $n\geq 1$, we have
\[
v_3(\phi_n(x,y)) \geq v_3(\phi_n)+n-1 \geq n-1-\frac{6n-4}{7} =
\frac{n-3}{7}>-1.
\]
Therefore $v_3(\phi_n(x,y))\geq 0$, which means that $\phi_n(x,y)$
can be evaluated in $\fg$ for all $n\geq 1$, and, in addition,
$v_3(\phi_n(x,y))\to\infty$ as $n\to\infty$ (independently of
$x,y$), which implies that the series $\phi(x,y)$ converges
uniformly in $\fg$ for $x,y\in\fg$. Similarly, the series
$\psi(x,y)$ can be evaluated term-by-term and converges uniformly
for all $x,y\in\fg$.

\mbr

Finally, we consider the case $p\geq 5$. Here an additional small
trick is needed. Put $\cK=E=\bQ_p(\sqrt{p})$ and
$H(x,y)=\frac{1}{\sqrt{p}}\cdot
CH(\sqrt{p}x,\sqrt{p}y)\in\cK\langle\langle x,y\rangle\rangle$.
Then \eqref{e:val-CH} implies that $v_p(H_n)=n/2+v_p(CH_n)-1/2\geq
0$ for all $n\geq 1$, and we have $H_1(x,y)=x+y$, which shows that
Lemma \ref{l:phi-psi} applies to $H(x,y)$. Changing notation
again, we let $\phi'(x,y),\,\psi'(x,y)$ be the formal Lie series
satisfying the conclusion of the lemma. Thus
\[
\frac{1}{\sqrt{p}} \log\bigl(\exp(\sqrt{p}x)\exp(\sqrt{p}y)\bigr)
= e^{\ad\phi'(x,y)}(x) + e^{\ad\psi'(x,y)}(y),
\]
which after a change of variables $z=\sqrt{p}x$, $w=\sqrt{p}y$ can
be rewritten as
\[
\log\bigl(\exp(z)\exp(w)\bigr) = e^{\ad\phi(z,w)}(z) +
e^{\ad\psi(z,w)}(w),
\]
where we have put
$\phi(z,w)=\phi'\bigl(\frac{z}{\sqrt{p}},\frac{w}{\sqrt{p}}\bigr)$
and
$\psi(z,w)=\psi'\bigl(\frac{z}{\sqrt{p}},\frac{w}{\sqrt{p}}\bigr)$.
The problem is that the coefficients of $\phi$ and $\psi$ lie a
priori in $\cK=\bQ_p(\sqrt{p})$. However, this is easy to fix as
follows. Let us introduce a $\bZ/2\bZ$-grading on
$\cK\langle\langle x,y\rangle\rangle$ by assigning degree $0$ to
every element of $\bQ_p$, and assigning degree $1$ to $x$, $y$ and
$\sqrt{p}$. With this convention, it is clear that $H(x,y)$ is
\emph{purely odd} (i.e., each $H_n(x,y)$ is odd). By looking at
the proof of Lemma \ref{l:phi-psi}, it is easy to see that the
formal Lie series $\phi'$ and $\psi'$ can be chosen to be purely
even. This implies that $\phi(z,w)$ and $\psi(z,w)$ have
coefficients in $\bQ_p$.

\mbr

The rest of the proof is the same as before. For all $z,w\in\fg$
and all $n\geq 1$, we have
\begin{eqnarray*}
v_p(\phi_n(z,w)) &\geq& v_p(\phi_n) + (n-1) = v_p(\phi'_n)
-\frac{n}{2}+(n-1) \\
 &\geq& -\frac{n-1}{p-1} - \frac{n}{2} + n-1 \geq
 \frac{n-3}{4}>-1,
\end{eqnarray*}
where we have used\footnote{Note that this argument would fail if
$p=3$: this is why we need Lemma \ref{l:ugly}.} the assumption
$p\geq 5$. Thus $\phi(z,w)$ can be evaluated term-by-term in $\fg$
and converges uniformly for all $z,w\in\fg$, and similarly for
$\psi(z,w)$.
\end{proof}

\mbr

\begin{thm}\label{t:orbmethod-not-2}
Assume that $p\geq 3$, let $G$ be either a $p$-group of nilpotence
class $<p$ or a uniform pro-$p$-group, and let $\fg=\Lie(G)$. Then
there exists a bijection $\Om\longleftrightarrow\chi_\Om$ between
$G$-orbits $\Om\subset\fg^*$ and characters of representations
$\rho\in\Gh$ such that Kirillov's character formula holds:
\begin{equation}\label{e:kir}
\chi_\Om(e^x) = \abs{\Om}^{-1/2}\cdot \sum_{f\in\Om} f(x) \qquad
\forall\,x\in\fg.
\end{equation}
\end{thm}
\begin{proof}
We will show that hypothesis (ii) of Theorem \ref{t:abstract}
holds for $\exp:\fg\rar{}G$. Let $\phi(x,y)$ and $\psi(x,y)$ be
the formal Lie series satisfying the conclusion of Proposition
\ref{p:phi-psi}. Then we obtain a continuous map of $\fg\times\fg$
to itself given by
$(x,y)\longmapsto(\xt,\yt)=\bigl(e^{\ad\phi(x,y)}(x),e^{\ad\psi(x,y)}(y)\bigr)$.
(This is a slight abuse of notation since $\xt$ depends on both
$x$ and $y$, and so does $\yt$.) This map satisfies the properties
mentioned in \S\ref{ss:formal-Lie}: on the one hand, $\xt$ and
$\yt$ are conjugate to $x$ and $y$, respectively, and on the other
hand, we have $CH(x,y)=\xt+\yt$ for all $x,y\in\fg$. We will now
use this information to show that
$\exp^*(f_1*f_2)=\exp^*(f_1)*\exp^*(f_2)$ for all
$f_1,f_2\in\Fun(G)^G$.

\mbr

Let us first assume that $G$ is finite. Then, due to the
$G$-invariance of $f_1$, $f_2$, we have
\begin{equation}\label{e:convol-G}
\bigl(\exp^*(f_1*f_2)\bigr)(z)
 \overset{\text{def}}{=} \frac{1}{\abs{\fg}}
\sum_{x,y\in\fg \ :\ \xt+\yt=z} f_1(e^x) f_2(e^y) =
\frac{1}{\abs{\fg}} \sum_{x,y\in\fg \ :\ \xt+\yt=z} f_1(e^{\xt})
f_2(e^{\yt})
\end{equation}
for all $z\in\fg$. On the other hand,
\begin{equation}\label{e:convol-g}
\bigl(\exp^*(f_1)*\exp^*(f_2)\bigr)(z) = \frac{1}{\abs{\fg}}
\sum_{x,y\in\fg \ :\ x+y=z} f_1(e^x) f_2(e^y).
\end{equation}
Thus it only remains to show that the map
$(x,y)\longmapsto(\xt,\yt)$ is a bijection of $\fg\times\fg$ onto
itself. However, this map is of the form
$(x,y)\longmapsto\bigl(x+A(x,y),y+B(x,y)\bigr)$, where $A$ and $B$
are Lie polynomials whose homogeneous components have degrees
$\geq 2$. Using induction on the nilpotence class of $\fg$, it is
easy to check that this map is injective (for the induction step,
let $\fz$ be the center of $\fg$ and note that the map descends to
a map of $(\fg/\fz)\times(\fg/\fz)$ to itself). Therefore it is
bijective because $\fg$ is finite.

\mbr

If $G$ is uniform the argument is similar. We only need to recall
that by the definition of $\Fun(G)$, there exists $r\in\bN$ such
that $f_1$ and $f_2$ are bi-invariant with respect to the open
subgroup $G^{p^r}=\exp(p^r\fg)$. Moreover, $\exp$ descends to a
bijection of $\fg/p^r\fg$ onto $G/G^{p^r}$, and the map
$(x,y)\longmapsto(\xt,\yt)$ descends to a bijection of
$(\fg/p^r\fg)\times(\fg/p^r\fg)$ onto itself. Thus equations
\eqref{e:convol-G} and \eqref{e:convol-g} remain valid with $\fg$
replaced by $\fg/p^r\fg$, and we see that
$\exp^*(f_1*f_2)=\exp^*(f_1)*\exp^*(f_2)$, as desired.
\end{proof}

\subsection{Uniform pro-$2$-groups}\label{ss:pro-2-groups} In this
subsection we assume that $p=2$, fix a uniform pro-$2$-group $G$,
and put $\fg=\Lie(G)$. As we have already mentioned, the
conclusion of Theorem \ref{t:orbmethod-not-2} may fail in this
case; indeed, even for an orbit $\Om\subset\fg^*$ of size $1$
formula \eqref{e:kir} may fail to define an irreducible character
of $G$ (cf.~\cite{jaikin-zapirain-zeta}). In view of Theorem
\ref{t:abstract}, this means that the pullback map
$\exp^*:\Fun(G)^G\rar{}\Fun(\fg^*)^G$ may not commute with
convolution. However, we do have a weaker positive result:
\begin{prop}\label{p:convol-2}
\begin{enumerate}[$($a$)$]
\item Assume that $[\fg,\fg]\subseteq 8\cdot\fg$. Given
$f_1,f_2\in\Fun(G)^G$, we have
$\exp^*(f_1*f_2)=\exp^*(f_1)*\exp^*(f_2)$ provided either $f_1$ or
$f_2$ is supported on $G^2$.
\item In general, we have
$\exp^*(f_1*f_2)=\exp^*(f_1)*\exp^*(f_2)$ for all
$f_1,f_2\in\Fun(G^2)^G$.
\end{enumerate}
\end{prop}
\begin{proof}
(a) The argument is rather similar to the one used in the previous
subsection. Consider the formal Lie series
$H(x,y)=\frac{1}{2}\cdot CH(2x,2y)$. We have $H_1(x,y)=x+y$, and
it follows from \eqref{e:val-CH} that $v_2(H_n)\geq 0$ for all
$n\geq 1$. Thus we can apply Lemma \ref{l:phi-psi} with $\cK=\bQ$,
and it yields formal Lie series $\phi'(x,y),\,\psi'(x,y)$
satisfying
\[
\frac{1}{2}\cdot\log\bigl(\exp(2x)\exp(2y)\bigr) =
e^{\ad\phi'(x,y)}(x)+e^{\ad\psi(x',y')}(y).
\]
We make the change of variables $z=2x$, $w=2y$ and rewrite the
last equation as
\[
\log\bigl(\exp(z)\exp(w)\bigr) = e^{\ad\phi(z,w)}(z) +
e^{\ad\psi(z,w)}(w),
\]
where $\phi(z,w)=\phi'\bigl(\frac{z}{2},\frac{w}{2}\bigr)$ and
$\psi(z,w)=\psi'\bigl(\frac{z}{2},\frac{w}{2}\bigr)$. Now if
$z,w\in\fg$, then for all $n\geq 1$,
\begin{equation}\label{e:val-2}
\begin{split}
v_2(\phi_n(z,w))&\geq v_2(\phi_n)+3(n-1)=v_2(\phi'_n)-n+3(n-1) \\
&\geq -(n-1)-n+3(n-1) = n-2.
\end{split}
\end{equation}
Here we have used the assumption $[\fg,\fg]\subseteq 8\cdot\fg$.
Similarly, $v_2(\psi_n(z,w))\geq n-2$ for all $n\geq 1$. This
means that the series $\sum_{n\geq 2}\phi_n(z,w)$ and $\sum_{n\geq
2}\psi_n(z,w)$ can be evaluated term-by-term and converge
uniformly for $z,w\in\fg$. Unfortunately, we cannot make sure that
both $\phi_1(z,w)$ and $\psi_1(z,w)$ are defined in $\fg$, because
by definition we must have
$[\phi_1(z,w),z]+[\psi_1(z,w),w]=\frac{1}{2}[z,w]$. However, in
view of the inductive construction of the series $\phi'$ and
$\psi'$ used in the proof of Lemma \ref{l:phi-psi}, we may assume
that, say, $\phi_1(z,w)=0$ and $\psi_1(z,w)=z/2$. This implies
that $\phi(z,w)$ is defined and converges uniformly in $\fg$ for
all $z,w\in\fg$, while $\psi(z,w)$ is defined and converges
uniformly in $\fg$ for $z\in 2\fg$ and $w\in\fg$.

\mbr

The rest of the proof is as before. Put
$(\xt,\yt)=\bigl(e^{\ad\phi(x,y)}(x),e^{\ad\psi(x,y)}(y)\bigr)$.
Then $(x,y)\longmapsto(\xt,\yt)$ is a map from $(2\fg)\times\fg$
to itself, and the argument used in the proof of Theorem
\ref{t:orbmethod-not-2} implies that
$\exp^*(f_1*f_2)=\exp^*(f_1)*\exp^*(f_2)$ if $f_1,f_2\in\Fun(G)^G$
and $f_1$ is supported on $G^2$. Since convolution of
$G$-invariant functions is commutative, the same formula holds if
instead $f_2$ is supported on $G^2$, completing the proof of (a).

\mbr

The proof of (b) is almost identical, except that \eqref{e:val-2}
has to be replaced by the following estimate, which is valid
whenever $[\fg,\fg]\subseteq 4\fg$ and $z,w\in 2\fg$:
\begin{equation*}
\begin{split}
v_2(\phi_n(z,w))&\geq v_2(\phi_n)+2(n-1)+n=v_2(\phi'_n)-n+2(n-1)+n \\
&\geq -(n-1)+2(n-1) = n-1\geq 0.
\end{split}
\end{equation*}
This means that $\phi(z,w)$, and similarly $\psi(z,w)$, can be
evaluated in $\fg$ term-by-term for all $z,w\in 2\fg$, and
converges uniformly for these values of $z,w$.
\end{proof}

We can now prove a version of the orbit method for uniform
pro-$2$-groups which is weaker than Theorem
\ref{t:orbmethod-not-2}, but suffices for some applications (see
Section \ref{s:brown}).
\begin{thm}\label{t:orbmethod-2}
Let $G$ be a uniform pro-$2$-group and $\fg=\Lie(G)$. For every
$G$-orbit $\Om\subset(2\fg)^*$, let $e'_\Om\in\Fun(2\fg)^G$ denote
the inverse Fourier transform of the characteristic function of
$\Om$, put $e_\Om=\log^*(e'_\Om)\in\Fun(G^2)^G$, and define
$\Gh_\Om\subset\Gh$ to be the collection of those $\rho\in\Gh$ on
which $e_\Om$ acts nontrivially\footnote{This means that the
linear operator $\rho(e_\Om):=\int_G \overline{e_\Om(g)}\rho(g)
d\mu_G(g)$ is nonzero.}. Then the following statements hold:
\begin{enumerate}[$($a$)$]
\item each $\Gh_\Om$ is finite, and $\Gh$ is the disjoint union of
the subsets $\Gh_\Om$;
\item if $\rho\in\Gh_\Om$ and $\chi_\rho$ is its character, then
$\chi_\rho\bigl\lvert_{G^2}$ is a multiple of $e_\Om$.
\end{enumerate}
\end{thm}
\begin{proof}
We use a modification of the argument that appeared in the proof
of Theorem \ref{t:abstract}. By construction, $e'_\Om$ is an
indecomposable idempotent in the algebra $\Fun(2\fg)^G$ (with
respect to the convolution defined using addition in $\fg$), and
Proposition \ref{p:convol-2}(b) implies that $e_\Om$ is an
indecomposable idempotent in $\Fun(G^2)^G$. Now we can think of
$\Fun(G^2)^G$ as a subalgebra of $\Fun(G^2)^{G^2}$ in the obvious
way, as well as a subalgebra of $\Fun(G)^G$ using extension by
zero. Therefore we can write
\[
e_\Om = \sum_{i=1}^m e_i = \sum_{j=1}^n f_j,
\]
where the $e_i$'s are indecomposable idempotents in
$\Fun(G^2)^{G^2}$ and the $f_j$'s are indecomposable idempotents
in $\Fun(G)^G$. By the proof of Theorem \ref{t:abstract}, each
$e_i$ (resp., $f_j$) corresponds to some $\pi_i\in\widehat{G^2}$
(resp., $\rho_j\in\Gh$) whose character is a multiple of $e_i$
(resp., $f_j$). It is clear that if $\rho\in\Gh$, then
$\rho(e_\Om)\neq 0$ if and only if $\rho\cong\rho_j$ for some $j$,
which implies that $\Gh_\Om=\{\rho_1,\rho_2,\dotsc,\rho_n\}$ is
finite, proving the first half of (a).

\mbr

Next let $\chi_i\in\Fun(G)^G$ be the character of the induced
representation $\eta_i:=\Ind_{G^2}^G\pi_i$. Since $G^2$ is normal
in $G$, it follows that $\chi_i$ is supported on $G^2$, so that we
can think of it as an element of $\Fun(G^2)^G$, and, moreover,
$\chi_i$ is a positive integral multiple of the sum of elements in
the orbit of $e_i$ under the $G$-conjugation action. In
particular, $\chi_i*e_\Om\neq 0$, and therefore
$\chi*e_\Om=\la_i\cdot e_\Om$ for some $\la_i\in\cst$, because
$e_\Om$ is a indecomposable idempotent in $\Fun(G^2)^G$. Hence we
must have $\chi_i*f_j=\la_i\cdot f_j$ for every $j$. Therefore the
$\rho_j$'s are precisely the irreducible constituents of $\eta_i$.
Now the Frobenius reciprocity implies that for each $1\leq j\leq
n$, the $\pi_i$'s are precisely the irreducible constituents of
$\rho_j\bigl\lvert_{G^2}$, which proves part (b).

\mbr

Finally, to finish the proof of (a), let $\rho\in\Gh$ be arbitrary, and let
$f\in\Fun(G)^G$ denote the corresponding indecomposable idempotent. There
exists a normal open subgroup $K\subset G$ such that $K\subseteq G^2$ and $f$
is bi-invariant with respect to $K$. Therefore $f$ is the pullback of an
indecomposable idempotent $\overline{f}$ of $\Fun(G/K)^{G/K}$. However, the
natural inclusion $\Fun(G^2/K)^{G/K}\hookrightarrow\Fun(G/K)^{G/K}$ is a
homomorphism of \emph{unital} algebras, which implies that $\overline{f}$ is a
summand of an indecomposable idempotent $\overline{e}$ of $\Fun(G^2/K)^{G/K}$.
Let $e\in\Fun(G^2)^G$ be the pullback of $\overline{e}$; it follows from
Proposition \ref{p:convol-2}(b) that $e=e_\Om$ for some $G$-orbit
$\Om\subset(2\fg)^*$, and the proof is complete.
\end{proof}

\subsection{Concluding remarks}\label{ss:conclusion} The orbit
method for uniform pro-$p$-groups was first studied by Roger Howe
\cite{howe}; he used the classical approach based on the notion of
a polarization. However, he did not treat the case $p=2$, and his
results for $p\geq 3$ are weaker than our Theorem
\ref{t:orbmethod-not-2} in that he has to impose an additional
requirement on $\fg$: namely, the Lie algebra $\widetilde{\fg}$
which has $\fg$ as the underlying $\bZ_p$-module and has the Lie
bracket defined by
$[x,y]_{\widetilde{\fg}}=\frac{1}{p}\cdot[x,y]_{\fg}$ must be
\emph{pro-nilpotent} (equivalently,
$\widetilde{\fg}/(p\cdot\widetilde{\fg})$ must be a nilpotent Lie
algebra over $\bF_p$). Thus, for example, Howe's result does not
apply to groups such as the kernel of the reduction modulo $p$
homomorphism $GL_n(\bZ_p)\rar{} GL_n(\bF_p)$ for any $p\geq 3$,
whereas our results do apply to them.

\mbr

The problem with the classical approach is that not every
polarization of $\fg$ corresponds to a subgroup of $G$, and Howe
imposed his assumption precisely to deal with it. However, Andrei
Jaikin-Zapirain showed in \cite{jaikin-zapirain-zeta} that Howe's
assumption can be removed by proving the existence of
polarizations satisfying some stronger conditions which allow the
classical method to be used. His Theorem 2.9 is equivalent to our
Theorem \ref{t:orbmethod-not-2}. Moreover, he also obtained a
result in the case $p=2$ (Theorem 2.12 in \emph{op.~cit.}) which
is stronger than our Theorem \ref{t:orbmethod-2}. On the other
hand, our result is already sufficient for some applications, as
we demonstrate in the next section.


\section{A $p$-adic analogue of Brown's theorem}\label{s:brown}

\subsection{The setup}\label{ss:brown-setup}
We warn the reader that our notation here will differ from that of
the first two sections. Namely, throughout the rest of the paper
we let $G$ be a \emph{$p$-adic nilpotent Lie group}, and $\fg$ its
Lie algebra, which is a finite dimensional nilpotent Lie algebra
over $\bQ_p$. For our purposes one does not need to know the
general definition of a $p$-adic Lie group; it suffices to think
of $G$ as the underlying topological space of $\fg$ (where the
topology on $\fg$ is induced by the standard topology on $\bQ_p$)
equipped with the operation given by the Campbell-Hausdorff series
$CH(x,y)$. (Here there is no question of $CH$ being well defined
or convergent, because $\fg$ is a Lie algebra over a field of
characteristic zero, and is nilpotent.) Thus $G$ is a locally
compact totally disconnected topological group. Recall also that a
choice of a nontrivial continuous additive character
$\psi:\bQ_p\to\cst$ allows one to identify $\fg^*$ with
$\Hom_{\bQ_p}(\fg,\bQ_p)$. The set of coadjoint orbits $\fg^*/G$
is equipped with the quotient of the natural topology of $\fg^*$.

\mbr

A complex representation $(\pi,V)$ of $G$, where $V$ is a vector
space over $\bC$ and $\pi:G\rar{}GL(V)$ is a homomorphism, is said
to be \emph{algebraic} (or \emph{smooth}) if for each $v\in V$ the
stabilizer $G^v=\{g\in G\st \pi(g)v=v \}$ is an open subgroup of
$G$. Note that $V$ need not have finite dimension in this
definition, and in fact most irreducible algebraic representations
of $G$ are infinite-dimensional. The isomorphism class of
$(\pi,V)$ will be denoted by $[(\pi,V)]$, and we write $\Gh$ for
the set of isomorphism classes of irreducible algebraic
representations of $G$. It is equipped with the Fell topology,
whose definition is recalled in \S\ref{ss:fell} below. Calvin
Moore proved \cite{moore} that there is a natural bijection
between $\fg^*/G$ and $\Gh$ (see the introduction). The main
result of this section is

\begin{thm}\label{t:brown}
The orbit method bijection $\fg^*/G\rar{}\Gh$ is a homeomorphism.
\end{thm}

Note that the continuity of this bijection is not difficult to
check using an argument similar to the one for real Lie groups
(but see also \S\ref{ss:continuity}). On the other hand, it is
rather nontrivial to prove that the bijection is open, and our
argument is based on a fact (Lemma \ref{l:union}) which does not
have an obvious analogue over $\bR$.

\subsection{Fell topology}\label{ss:fell} We recall the definition given
in \cite{gelfand-kazhdan}. For an irreducible algebraic
representation $(\pi,V)$ of $G$, choose $n\in\bN$, vectors
$v_1,\dotsc,v_n\in V$, linear functionals $\xi_1,\dotsc,\xi_n\in
V^*$, a compact set $B\subset G$, and a real number $\eps>0$, and
define
\[
\sU(\pi,V,B,v_j,\xi_j,\eps) \subseteq\Gh
\]
to be the set of isomorphism classes $[(W,\rho)]\in\Gh$ such that
there exist $w_1,\dotsc,w_n\in W$ and $\eta_1,\dotsc,\eta_n\in
W^*$ with the property
\[
\Bigl\lvert \bigl\langle \xi_i,\pi(g)v_i \bigr\rangle -
\bigl\langle \eta_i, \rho(g) w_i \bigr\rangle \Bigr\rvert < \eps
\qquad \forall\,g\in B,\ 1\leq i\leq n.
\]
Sets of the form $\sU(\pi,V,B,v_j,\xi_j,\eps)$ are defined to be a
basis of neighborhoods of the point $[(\pi,V)]\in\Gh$, which
uniquely determines a topology on $\Gh$, called the \emph{Fell
topology}. To understand it, we begin with the following
\begin{lem}\label{l:union}
If $\fg$ is a finite dimensional nilpotent Lie algebra over
$\bQ_p$, then $\fg$ can be written as the union of an increasing
sequence of open uniform Lie subalgebras:
\begin{equation}\label{e:union}
\fk_1 \subseteq \fk_2 \subseteq \fk_3 \subseteq \dotsb
\subseteq\fg, \qquad \fg=\bigcup_{j\geq 1} \fk_j.
\end{equation}
\end{lem}
\begin{proof}
Let $x_1,\dotsc,x_N$ be a basis of $\fg$ over $\bQ_p$. For every
$j\in\bN$, consider the set of all elements of $\fg$ of the form
\[
[y_1,[y_2,[y_3,[\dotsb [y_{t-1},y_t]\dotsb ]]]],
\]
where $t\geq 1$ is arbitrary and each $y_i$ is of the form $p^{-j}
x_k$ for some $1\leq k\leq N$. Since $\fg$ is nilpotent, only
finitely many of these iterated commutators are nonzero, and hence
their $\bZ_p$-span, call it $\fk'_j$, is a free $\bZ_p$-submodule
of $\fg$ of finite rank. Moreover, $\fk'_j$ is closed under the
Lie bracket by definition. By construction,
$\fk'_j\subseteq\fk'_{j+1}$ for all $j\geq 1$, and
$\fg=\fk'_1\cup\fk'_2\cup\dotsb$. Let $\fk_j=p\cdot\fk'_j$ (resp.,
$\fk_j=4\cdot\fk'_j$ if $p=2$); this is clearly a uniform Lie
algebra, and since $\fg$ is a vector space over $\bQ_p$, it
follows that \eqref{e:union} holds. Finally, each $\fk_j$ is open
because $p^{1-j}x_m\in\fk_j$ for all $1\leq m\leq N$.
\end{proof}

\subsection{Some notation}\label{ss:brown-notation} An obvious consequence of
Lemma \ref{l:union} is that one obtains the same topology on $\Gh$
by restricting the compact set $B$ in the definition of the Fell
topology to be an arbitrary open uniform subgroup $K$ of $G$. Now
either Theorem \ref{t:orbmethod-not-2} or Theorem
\ref{t:orbmethod-2} applies to $K$. In order to make the argument
below independent of $p$, let us define $\al=2$ if $p=2$ and
$\al=1$ if $p\geq 3$. Put $\fk=\Lie(K)\subset\fg$, and let
$\Om_0\subset(\al\fk)^*$ be a $K$-orbit. If $p=2$, then a finite
subset $\Kh_{\Om_0}\subset\Kh$ was defined in Theorem
\ref{t:orbmethod-2}. For $p\geq 3$, we let $\Kh_{\Om_0}\subset\Kh$
to be the singleton subset consisting of the irreducible
representation of $K$ corresponding to $\Om_0\subset\fk^*$. In
either of the two cases, we define $e_{\Om_0}\in\Fun(K^\al)^K$ be
the pullback via $\log:K^\al\rar{}\al\fk$ of the inverse Fourier
transform of the characteristic function of $\Om_0$, and Theorems
\ref{t:orbmethod-not-2} and \ref{t:orbmethod-2} imply that if
$\rho\in\Kh_{\Om_0}$ and $\chi_\rho$ is its character, then
$\chi_\rho\bigl\lvert_{K^\al}$ is a multiple $e_{\Om_0}$.

\mbr

As the last piece of notation, if $\pi$ is any complex continuous
representation of $K$, we will denote by $\Supp(\pi)\subseteq\Kh$
the collection of all irreducible constituents of $\pi$.

\subsection{Proof of the difficult part of Theorem \ref{t:brown}}\label{ss:brown-proof} Given $f\in\fg^*$,
we will prove that the orbit method bijection
$\fg^*/G\rarab{\simeq}\Gh$ is open at the point $\Om_f\in\fg^*/G$,
where $\Om_f$ denotes the $G$-orbit of $f$. Consider an open
neighborhood of $f$ in $\fg^*$. By shrinking it if necessary, we
may assume, thanks to Lemma \ref{l:union}, that it is of the form
\[
\sV(f,K)=\Bigl\{f'\in\fg^* \,\Bigl\lvert\,
f'\big\vert_{\al\fk}=f\big\vert_{\al\fk} \Bigr\},
\]
where $K\subset G$ is an open uniform subgroup with Lie algebra
$\fk$ and $\al$ is as in \S\ref{ss:brown-notation}. We assume from
now on that $K$ is fixed. It suffices to check that the image of
$\sV(f,K)$ under the orbit method map $\fg^*\rar{}\Gh$ contains an
open neighborhood of $[(\pi_f,V_f)]$ with respect to the Fell
topology, where $(\pi_f,V_f)$ denotes an irreducible algebraic
representation of $G$ corresponding to $\Om_f$. The proof rests on
the following
\begin{prop}\label{p:restriction}
Let $\Om\subset\fg^*$ be a $G$-orbit, let $\Om_0\subset(\al\fk)^*$
be a $K$-orbit, and let $\pi$ denote the irreducible algebraic
representation of $G$ corresponding to $\Om$. Then
$\Supp(\pi\bigl\lvert_K)\cap\Kh_{\Om_0}\neq\varnothing$ if and
only if $\Om_0$ is contained in the image of $\Om$ under the
restriction map $\res:\fg^*\to(\al\fk)^*$.
\end{prop}
\begin{proof}
We use \cite{gelfand-kazhdan}, \S1.2. Recall that the character
$c(\pi)$ of $\pi$ is not a function on $G$, but rather a
distribution, defined by
\begin{equation}\label{e:character}
\pair{c(\pi),t} = \tr\left[ \int_G \overline{t(g)}\pi(g)dg
\right], \qquad t\in C_0^\infty(G),
\end{equation}
where $C^\infty_0(G)$ is the space of locally constant functions $G\rar{}\bC$
with compact support, and $dg$ denotes a fixed Haar measure on $G$. (The
complex conjugation appears in the formula above for consistency with our orbit
method for uniform pro-$p$-groups.) Moreover, Kirillov's character formula in
this context implies that $\exp^*(c(\pi))$ is the inverse Fourier transform of
a suitably normalized $G$-invariant measure on $\fg^*$ supported on $\Om$. Let
$e_0=e_{\Om_0}\in\Fun(K^\al)^K\subset C_0^\infty(G)$ with the notation of
\S\ref{ss:brown-notation}. By definition,
$\Supp(\pi\bigl\lvert_K)\cap\Kh_{\Om_0}\neq\varnothing$ if and only if $e_0$
acts nontrivially on $\pi\bigl\lvert_K$, which in turn is equivalent to
$\pair{c(\pi),e_0}\neq 0$ because $e_0$ is an idempotent in $\Fun(K)$. Since
the Fourier transform is an isometry, we see that $\pair{c(\pi),e_0}\neq 0$ if
and only if $\res^{-1}(\Om_0)\cap\Om\neq\varnothing$, which proves the
proposition (since $\res(\Om)$ is obviously $K$-stable).
\end{proof}

\mbr

Using the notation preceding the statement of the proposition, let
$f_0=f\big\vert_{\al\fk}$, let $\Om_0\subset(\al\fk)^*$ be the
$K$-orbit of $f_0$, and let $e_0=e_{\Om_0}$. Then
$\Om_0\subseteq\res(\Om)$, so by (the proof of) Proposition
\ref{p:restriction}, we have $\pi_f(e_0)\neq 0$. In particular,
there exist $v\in V_f$ and $\xi\in V_f^*$ such that
$\pair{\xi,\pi_f(e_0)v}=1$. Define $\eps=\left( \int_K
\abs{e_0(k)} d\mu_K(k) \right)^{-1}$, where $\mu_K$ is the
standard Haar measure on $K$ of total mass $1$. It is clear that
the following result implies that the orbit method bijection
$\fg^*/G\rar{}\Gh$ is open:
\begin{prop}\label{p:open}
The open neighborhood $\sU(\pi_f,V_f,K,v,\xi,\eps)$ of
$[(\pi_f,V_f)]$ in $\Gh$ is contained in the image of $\sV(f,K)$
under the map $\fg^*\rar{}\Gh$.
\end{prop}
\begin{proof}
Suppose that $[(\rho,W)]\in\sU(\pi_f,V_f,K,v,\xi,\eps)$. By
definition, there exist $w\in W$ and $\eta\in W^*$ such that
$\bigl\lvert \pair{\eta,\rho(g)w} - \pair{\xi,\pi_f(g)v}
\bigr\rvert < \eps$ for all $g\in K$. Multiplying by the function
$\abs{e_0(g)}$, integrating with respect to $\mu_K$ and using the
definition of $\eps$, we obtain $\bigl\lvert
\pair{\eta,\rho(e_0)w} - \pair{\xi,\pi_f(e_0)v} \bigr\rvert < 1$.
By construction, this forces $\pair{\eta,\rho(e_0)w}\neq 0$,
whence $\rho(e_0)\neq 0$, i.e.,
$\supp(\rho\bigl\lvert_K)\cap\Kh_{\Om_0}\neq\varnothing$. By
Proposition \ref{p:restriction}, if $\Om'\subset\fg^*$ is the
$G$-orbit corresponding to $\rho$, then $\res(\Om')$ contains
$f_0$, i.e., $\Om'$ meets $\sV(f,K)$.
\end{proof}

\subsection{Proof of continuity (sketch)}\label{ss:continuity} We
conclude by briefly explaining how Lemma \ref{l:union} can be used
to prove that the orbit method bijection $\fg^*/G\rar{}\Gh$ is
continuous. We use the notation of \S\ref{ss:fell}. Fix
$(\pi,V)\in\Gh$, let $\Om\subset\fg^*$ be the corresponding
$G$-orbit, and consider a ``standard'' open neighborhood of
$[(\pi,V)]$ of the form $\sU=\sU(\pi,V,K,v_j,\xi_j,\eps)$, where
$K\subset G$ is a uniform pro-$p$-subgroup. In view of Lemma
\ref{l:union} it suffices to show that the inverse image of $\sU$
in $\fg^*/G$ contains a neighborhood of $\Om$. To this end, let
$W\subseteq V$ be a finite dimensional $K$-invariant subspace
containing all the $v_j$'s, let $\vartheta$ denote the (possibly
reducible) representation of $K$ afforded by $W$, and let
$\chi_1,\dotsc,\chi_r\in\Fun(K)^K\subset C_0^\infty(G)$ denote the
characters of the irreducible constituents of $\vartheta$. It is
clear that if $(\pi',V')\in\Gh$ is such that $\pi'\bigl\lvert_K$
contains a $K$-subrepresentation isomorphic to $\vartheta$, then
$[(\pi',V')]\in\sU$.

\mbr

Now fix a Haar measure $dg$ on $G$. For every $G$-orbit
$\Om'\subset\fg^*$, let $\pi_{\Om'}$ be the corresponding
representation of $G$, let $c(\pi_{\Om'})$ be the character of
$\pi_{\Om'}$ defined by \eqref{e:character}, and let $\mu_{\Om'}$
denote the $G$-invariant measure on $\fg^*$ supported on $\Om'$
whose inverse Fourier transform is equal to
$\exp^*(c(\pi_{\Om'}))$. Furthermore, let $\fk=\Lie(K)\subset\fg$,
let $\res:\fg^*\rar{}\fk^*$ denote the restriction map, and let
$\nu_i:\fk^*\to\bC$ denote the Fourier transform of
$\exp^*(\chi_i)\in\Fun(\fk)^K$. It is not hard to check that as
$\Om'\in\fg^*/G$ varies, the conditions
$\bigl\langle\nu_i,\res_*\mu_{\Om'}\bigr\rangle\geq
\bigl\langle\nu_i,\res_*\mu_{\Om}\bigr\rangle$, $1\leq i\leq r$,
define an open subset $\sV\subseteq\fg^*/G$. If $\Om'\in\sV$, then
applying the inverse Fourier transform and Kirillov's character
formula for the group $G$, we see that the multiplicity of each
$\chi_i$ in $\pi_{\Om'}\bigl\lvert_K$ is at least its multiplicity
in $\pi\bigl\lvert_K$, whence $\vartheta$ is isomorphic to a
subrepresentation of $\pi_{\Om'}\bigl\lvert_K$. By the previous
paragraph, $\sV$ is contained in the inverse image of $\sU$ in
$\fg^*/G$, completing the proof.

\end{document}